\documentclass[12pt]{report}
\usepackage{amssymb,latexsym,epsfig,mathrsfs,graphpap,graphics,amsmath,amssymb}
\usepackage{amstext}
\usepackage{amsthm}
\usepackage[utf8]{inputenc}
\usepackage[english]{babel}
\usepackage{helvet}
\usepackage{courier}

\theoremstyle{Definition}

\theoremstyle{remark}

\numberwithin{equation}{section}

\begin{document}

\begin{flushleft}
 {\bf\Large { Uncertainty Principles for  Quaternion\\ \parindent=0mm \vspace{.1in}Linear Canonical S-Transform}}

\parindent=0mm \vspace{.2in}

{\bf{M. Younus Bhat$^{1},$ and Aamir H. Dar$^{2}$ }}
\end{flushleft}

\parindent=0mm \vspace{.1in}
{{\it Department of  Mathematical Sciences,  Islamic University of Science and Technology Awantipora, Pulwama, Jammu and Kashmir 192122, India.\\
 $^{1}$ gyounusg@gmail.com and $^{2}$ ahdkul740@gmail.com}}

\parindent=0mm \vspace{.2in}
{\bf{Abstract:}}  In this paper, we introduce the notion of quaternion linear canonical S-transform(Q-LCST)
which is an extension of the linear canonical S-transform. Firstly, we study the fundamental properties of quaternion linear canonical S-transform(Q-LCST) and then establish some basic results including orthogonality relation and reconstruction formula. Finally, we derive the associated Heisenberg’s uncertainty inequality and the corresponding logarithmic version for quaternion linear canonical S-transform(Q-LCST).

\parindent=0mm \vspace{.2in}

\parindent=0mm \vspace{.2in}

\parindent=0mm \vspace{.2in}
{\bf{1. Introduction}}

\parindent=0mm \vspace{.2in}
The quaternion linear canonical transform (QLCT), as a generalized form of the quaternion Fourier transform, is a powerful analyzing tool in image and signal processing.
A nontrivial generalization of the classical real and complex Fourier transform is the quaternion Fourier transform. This generalization using quaternion algebra has been  of great interest to the mathematicians and physicists for so many years. Many of the properties of the Fourier transform hold good in the new setting however others needed to be modified. It was R. G. Stockwell who constructed the novel Stockwell transform by using wavelet transform and the classical windowed Fourier transform. For more details we refer to [1, 2, 3, 9] and references therein. 

\parindent=0mm \vspace{.2in}
The well known Heisenberg's  uncertainty principle was  first introduced by German physicist Heisenberg. This principle  plays a vital  role in signal and image processing, physics and mathematics and other allied subjects. To be precise, it states that  a signal or function  cannot be located in both time and frequency domains, simultaneously and sharply. while as on the other side, the linear canonical transform (LCT) is a linear integral transformation with three free parameters, and it has been widely used in physical optics and signal processing. Indeed, it is also regarded as a generalization of many mathematical transforms such as the Fourier transform, the fractional Fourier transform and others. For further details we suggest [4, 5, 6, 7, 8, 10].

\parindent=0mm \vspace{.2in}
Continuing the studies of the extensions of the linear canonical transform, we introduced the notion of quaternion linear canonical S-transform(Q-LCST)
which is an extension of the linear canonical S-transform. Firstly, we study the fundamental properties of quaternion linear canonical S-transform(Q-LCST) and then establish some basic results including orthogonality relation and reconstruction formula. Finally, we derive the associated Heisenberg’s uncertainty inequality and the corresponding logarithmic version for quaternion linear canonical S-transform(Q-LCST). The rest of the paper is organised as follows:

\parindent=0mm \vspace{.2in}
In Section 2, we discuss some preliminary results and definitions which are used in subsequent sections.  In Section 3, we introduce the well known properties of Q-LCST like orthogonality, reconstruction formula and many more. In Section 4, we establish the well known Heisenberg’s uncertainty inequality for quaternion linear canonical S-transform. 

\parindent=0mm \vspace{.2in}
{\bf{Methods/Experiments}}

\parindent=0mm \vspace{.2in}
{\bf{2. Preliminaries  }}

\parindent=0mm \vspace{.2in}
In this section we review Linear canonical transform(LCT) , Linear canonical S-transform(LCST) and Quaternion linear canonical transform(QLCT).

\parindent=0mm \vspace{.1in}
{\it { 2.1.Linear canonical transform}}

\parindent=0mm \vspace{.1in}

 Let $\mu_1,\mu_2,\mu_3$ denote the three imaginary units in the quaternion algebra\cite{26OLCT} then the linear canonical transform (LCT) of any signal $f\in L^2(\mathbb R)$ is defined as
 \begin{align*}
\mathcal L_M[f](u)=\int_{\mathbb R}f(x)K^{\mu_1}_M(x,u)dx\tag{2.1}
\end{align*}
and its inversion is given by
\begin{align*}
f(x)&=\mathcal {L}^{-1}_M\{\mathcal L_M[f]\}(x)\\
&=\int_{\mathbb R}\mathcal L_M[f](u)K^{-\mu_1}_{M}(x,u)du\tag{2.2}
 \end{align*}
where $u,x\in \mathbb R,M=(A,B,C,D)$ is a uni-modular matrix parameter, and $K^{\mu_1}_M(x,u)$ is so-called the kernel of the LCT and is given by
\begin{align*}
K^{\mu_1}_M(x,u)=\begin{cases}
		\frac{1}{\sqrt{2\pi B}}e^{\frac{\mu_1}{2}\left(\frac{A}{B}x^2-\frac{2}{B}xu+\frac{D}{B}u^2-\frac{\pi}{2} \right) },   &B\neq0  \\
		\sqrt{D}e^{\mu_1\frac{CD}{2}u^2}\delta(x-Du),    &B=0
		\end{cases}\tag{2.3}
\end{align*}

\parindent=0mm \vspace{.1in}

{\it { 2.2.Linear Canonical S-Transform }}

\parindent=0mm \vspace{.1in}
The Linear Canonical S-Transform(ST) is known as a hybrid form of the S-transform. This transform is an extension of the LCT and ST,which gives a time-linear canonical domain joint expression. It is defined as:

For any non-zero window function $\psi\in L^2(\mathbb R)$ the LCST of any signal $f\in L^2(\mathbb R)$ with respect to $\psi$ is defined by
 \begin{align*}
\mathcal S_\psi^M[f](u,w)&=\mathcal {L}_M\left\{f(x)\overline{\psi(u-x,w)}\right\}(w,u)\\
&=\int_{\mathbb R}f(x)\overline{\psi(u-x,w)}K^{\mu_1}_{M}(x,u)dx\tag{2.4}
 \end{align*}

 and its inverse is given by

 \begin{align*}
 f(x)\overline{\psi(u-x,w)}&=\mathcal {L}^{-1}_M\{\mathcal S_\psi^Mf(u,w)\}\\
&=\int_{\mathbb R}\mathcal  S_\psi^Mf(u,w)K^{-\mu_1}_{M}(x,u)dw\tag{2.5}
 \end{align*}

 where $K^{\mu_1}_M(x,u)$ is given by (2.3 ) for $B\ne 0.$

 \parindent=0mm \vspace{.1in}

 It is worth noting that when $M=(0,1,-1,0),$ the LCST boils down S-transform.

 \parindent=0mm \vspace{.1in}

 {\it {2.3.Quaternion linear canonical transform}}

\parindent=0mm \vspace{.1in}
The QLCT is a generalization of the
LCT in the frame of quaternion algebra.The QLCT of a signal $f\in L^2(\mathbb R^2,\mathbb H)$ is defined as
\begin{align}
	\mathcal{L}_{M_1,M_2}^{\mathbb{H}}\{f\}(\mathbf{u})=
\int_{\mathbb{R}^2}K_{M_1}^{\mu_1}(x_1,u_1)f(\mathbf{x})K_{M_2}^{\mu_2}(x_2,u_2)d\mathbf{x},\tag{2.6}
	\end{align}	
and its inverse is given by
\begin{align*}
f(\mathbf{x})&=\mathcal {L}^{-1}_{M_1,M_2}\{\mathcal L^{\mathbb H}_{M_1,M_2}[f]\}(\mathbf {x})\\
&=\int_{\mathbb R^2}K^{-\mu_1}_{M_1}(u_1,x_1)\mathcal L^{\mathbb{H}}_{M_1,M_2}[f](\mathbf{u})K^{-\mu_2}_{M_2}(u_2,x_2)d\mathbf{u}\tag{2.7}
 \end{align*}
	where $\mathbf{u}=(u_1,u_2)$,$\mathbf{x}=(x_1,x_2)$ ,$M_s=(A_s,B_s,C_s,D_s)$ be a matrix parameter satisfying $det (M_s ) = 1$ and  the kernel
	\begin{align}
		\begin{split}
			 K_{M_s}^{\mu_s}(x_s,u_s)=\begin{cases}
		\frac{1}{\sqrt{2\pi B_s}}e^{\mu_s\left(\frac{A_s}{2B_s}x_s^2-\frac{x_su_s}{B_s}+\frac{D_s}{2B_s}u_s^2-\frac{\pi}{4} \right) },   &B_s\neq0  \\
		\sqrt{D_s}e^{\mu_s\frac{C_sD_s}{2}u_s^2}\delta(x_{s}-D_{s}u_{s}),    &B_s=0
		\end{cases}
		\end{split}\tag{2.8}
	\end{align}
for $s=1,2,\delta(x)$ represents the Dirac function.	

 \parindent=0mm \vspace{.2in}
{\bf{Results and Discussion }}

 \parindent=0mm \vspace{.2in}
{\bf{3.The Quaternion Linear Canonical S-Transform(Q-LCST)   }}

 \parindent=0mm \vspace{.1in}

 The linear canonical S-transform (LCST) is a linear integral
transformation having more than one free parameters and includes
the traditional Stockwell spectrum and the fractional Stockwell spectrum as its special cases.LCST is widely used in various fields viz: non-stationary signal representation, detection, parameters estimation,
 filter design and etc.This section will consider generalizing the
LCST using the quaternion algebra.This extension is then
called the quaternion linear canonical S-transform which is denoted by Q-LCST.Moreover, several
basic properties of them are investigated.

\parindent=0mm \vspace{.1in}

 {\it {3.1.Definition of Quaternion linear canonical S-transform}}

 \parindent=0mm \vspace{.1in}

 Based on the definition of the linear canonical S-transform(LCST), we obtain a definition of the quaternion linear canonical S-transform(Q-LCST) by replacing
the kernel of the linear canonical transform(LCT) with the kernel of the quaternion linear canonical transform(QLCT) in the classical definition of linear canonical S-transform(LCST) as follows.

\parindent=0mm \vspace{.1in}

{\bf Definition 3.1:} Let $\Psi\in L^2(\mathbb R^2,\mathbb H)$ be a non zero quaternion window function satisfying
 $$\int_{\mathbb R^2}\Psi(\mathbf{x})d\mathbf{x}=1\eqno(3.1)$$
then 2D quaternion linear canonical S-transform(Q-LCST)of quaternion signal $f\in L^2(\mathbb R^2,\mathbb H)$ with respect $\Psi$ is defined by
$$\mathbb S^{\mathbb H}_{\Psi,M_1,M_2}[f](\mathbf{u},\mathbf{w})=\int_{\mathbb R^2}K_{M_1}^{\mu_1}(x_1,w_1)f(\mathbf{x})\overline{\Psi(\mathbf{u}-\mathbf{x},\mathbf{w})}K_{M_2}^{\mu_2}(x_2,w_2)d\mathbf{x}\eqno(3.2)$$

where $\mathbf{u}=(u_1,u_2)\in{\mathbb R^2},\mathbf{w}=(w_1,w_2)\in{\mathbb R^2}$,$\mathbf{x}=(x_1,x_2)$ ,$M_s=(A_s,B_s,C_s,D_s)\in{\mathbb R^2}$ be a matrix parameter satisfying $det (M_s ) = 1$ and  the kernel
	\begin{align}
		\begin{split}
			 K_{M_s}^{\mu_s}(x_s,w_s)=\begin{cases}
		\frac{1}{\sqrt{2\pi B_s}}e^{\mu_s\left(\frac{A_s}{2B_s}x_s^2-\frac{x_sw_s}{B_s}+\frac{D_s}{2B_s}w_s^2-\frac{\pi}{2} \right) },   &B_s\neq0  \\
		\sqrt{D_s}e^{\mu_s\frac{C_sD_s}{2}w_s^2}\delta(x_{s}-D_{s}w_{s}),    &B_s=0
		\end{cases}
		\end{split}\tag{3.3}
	\end{align}
for $s=1,2,\delta(x)$ represents the Dirac function.

\parindent=0mm \vspace{.2in}

Note that the Q-LCST the case when $B_s = 0,s=1,2$
is not interesting because it is essentially a multiplication by a
quaternion chirp.Hence, without loss of generality,
we set $B_s\ne0,s=1,2$ in this paper.

It is worth noting that the Q-LCST  gives birth to some new time-frequency transforms
which are not yet reported in the open literature:

\begin{itemize}
\item For $M_s=(cos\theta,sin\theta,-sin\theta,cos\theta),\theta\ne n\pi$ we can obtain a novel quaternion fractional S-transform.
\item For $M_s=(1,B,0,1),B\ne 0$ we can obtain a new transform namely the Quaternion Fresnel S-transform.

\item For $M_s=(0,1,-1,0)$ the Q-LCST boils down to Quaternion S-transform.
\end{itemize}

\parindent=0mm \vspace{.1in}

From definition 3.1 we have

$$\mathbb S^{\mathbb H}_{\Psi,M_1,M_2}[f](\mathbf{u},\mathbf{w})=\mathcal L^{\mathbb H}_{M_1,M_2}\left\{f(\mathbf{x})\overline{\Psi(\mathbf{u}-\mathbf{x},\mathbf{w})}\right\}(\mathbf{w},\mathbf{u})\eqno(3.4)$$
applying inverse QLCT to(3.4),we have
\begin{align*}
f(\mathbf{x})\overline{\Psi(\mathbf{u}-\mathbf{x},\mathbf{w})}&=\mathcal L^{-1}_{M_1,M_2}\left\{\mathbb S^{\mathbb H}_{\Psi,M_1,M_2}[f](\mathbf{u},\mathbf{w})\right\}\\
&=\int_{\mathbb R^2}K_{M_1}^{-\mu_1}(x_1,w_1)\mathbb S^{\mathbb H}_{\Psi,M_1,M_2}[f](\mathbf{u},\mathbf{w})K_{M_2}^{-\mu_2}(x_2,w_2)d\mathbf{w}\tag{3.5}\\
\end{align*}

In the sequel, our intention is to study the fundamental properties of the proposed Q-LCST(3.2).

\parindent=0mm \vspace{.1in}

 {\it {3.2.Some Properties of Q-LCST}}

 \parindent=0mm \vspace{.1in}

 In this subsection, we discuss several basic properties of the Q-LCST. These
properties play important roles in signal representation.Most of the properties of the quaternion S-transform can be established in the Q-LCST domain with some modifications in quaternion windowed linear canonical transform.However, it is clearly visible that the properties of the
Q-LCST like shift,modulation and etc. vary from the corresponding properties of the quaternion windowed linear
canonical transform.

\parindent=0mm \vspace{.1in}

{\bf Property 1.} {\it (Multiplication operator) $ \left[\mathbb S^{\mathbb H}_{\Psi,M_1,M_2}\right]$ is a multiplication operator if  $\Psi(\mathbf{u-x},\mathbf{w})=\Gamma(\mathbf{w})$ i.e independent of $\mathbf{x}$. }

\parindent=0mm \vspace{.1in}

{\bf Proof.} From (3.2), we have

$$\mathbb S^{\mathbb H}_{\Psi,M_1,M_2}[f](\mathbf{u},\mathbf{w})=\int_{\mathbb R^2}K_{M_1}^{\mu_1}(x_1,w_1)f(\mathbf{x})\overline{\Psi(\mathbf{u}-\mathbf{x},\mathbf{w})}K_{M_2}^{\mu_2}(x_2,w_2)d\mathbf{x}\eqno(3.6)$$

Substituting $\Psi(\mathbf{u-x},\mathbf{w})=\Gamma(\mathbf{w})$ in (3.6),we obtain
\begin{align*}
\mathbb S^{\mathbb H}_{\Psi,M_1,M_2}[f](\mathbf{u},\mathbf{w})&=\int_{\mathbb R^2}K_{M_1}^{\mu_1}(x_1,w_1)f(\mathbf{x})\Gamma(\mathbf{w})K_{M_2}^{\mu_2}(x_2,w_2)d\mathbf{x}\\\\
&=\Gamma(\mathbf{w})\int_{\mathbb R^2}K_{M_1}^{\mu_1}(x_1,w_1)f(\mathbf{x})K_{M_2}^{\mu_2}(x_2,w_2)d\mathbf{x}\\\\
&=\Gamma(\mathbf{w})\mathcal L^{\mathbb H}_{M_1,M_2}[f](\mathbf{w}).\tag{3.7}
\end{align*}
So  $\left[\mathbb S^{\mathbb H}_{\Psi,M_1,M_2}\right]$ is a multiplication operator.

\parindent=0mm \vspace{.1in}

{\bf Example 1.} If $\Psi(\mathbf{u-x},\mathbf{w})=\Gamma(\mathbf{w})=1$ then (3.7) reduces to QLCT $$\mathbb S^{\mathbb H}_{\Psi,M_1,M_2}[f](\mathbf{u},\mathbf{w})=\mathcal L^{\mathbb H}_{M_1,M_2}[f](\mathbf{w}).\eqno(3.8)$$

\parindent=0mm \vspace{.1in}

{\bf Example 2.} If $\Psi(\mathbf{u-x},\mathbf{w})=\Gamma(\mathbf{x})$ i.e $\Psi$ is dependent on $\mathbf{x}$ then  from (3.6), we have
\begin{align*}
\mathbb S^{\mathbb H}_{\Psi,M_1,M_2}[f](\mathbf{u},\mathbf{w})&=\int_{\mathbb R^2}K_{M_1}^{\mu_1}(x_1,w_1)f(\mathbf{x})\Gamma(\mathbf{x})K_{M_2}^{\mu_2}(x_2,w_2)d\mathbf{x}\\\\
&=\int_{\mathbb R^2}K_{M_1}^{\mu_1}(x_1,w_1)g(\mathbf{x})K_{M_2}^{\mu_2}(x_2,w_2)d\mathbf{x}\\\\
&=\mathcal L^{\mathbb H}_{M_1,M_2}[g](\mathbf{w}).\tag{3.9}
\end{align*}
Where $g(\mathbf{x})=f(\mathbf{x})\Gamma(\mathbf{t})$
thus we see Q-LCST  boils down to QLCT.\\\\

\parindent=0mm \vspace{.4in}

{\bf Property 2.} {\it If $\Psi\in L^2(\mathbb R^2,\mathbb H)$ such that $\int_{\mathbb R^2}\Psi(\mathbf{u-x},\mathbf{w})d\mathbf{u}=1$ then $$\int_{\mathbb R^2}\mathbb S^{\mathbb H}_{\Psi,M_1,M_2}f(\mathbf{u,w})d\mathbf{u}=\mathcal L^{\mathbb H}_{M_1,M_2}[f](\mathbf{w})\eqno(3.10)$$}

\parindent=0mm \vspace{.1in}

{\bf Proof.} From (3.2), we have
\begin{align*}
\int_{\mathbb R^2}\mathbb S^{\mathbb H}_{\Psi,M_1,M_2}[f](\mathbf{u},\mathbf{w})&=\int_{\mathbb R^2}\int_{\mathbb R^2}K_{M_1}^{\mu_1}(x_1,w_1)f(\mathbf{x})\overline{\Psi(\mathbf{u}-\mathbf{x},\mathbf{w})}K_{M_2}^{\mu_2}(x_2,w_2)d\mathbf{x}d\mathbf{x}\\
&=\int_{\mathbb R^2}K_{M_1}^{\mu_1}(x_1,w_1)f(\mathbf{x})K_{M_2}^{\mu_2}(x_2,w_2)\int_{\mathbb R^2}\overline{\Psi(\mathbf{u}-\mathbf{x},\mathbf{w})}d\mathbf{u}d\mathbf{x}\\
&=\int_{\mathbb R^2}K_{M_1}^{\mu_1}(x_1,w_1)f(\mathbf{x})K_{M_2}^{\mu_2}(x_2,w_2)d\mathbf{x}\\
&=\mathcal L^{\mathbb H}_{M_1,M_2}[f](\mathbf{w})
\end{align*}
which completes proof.\qquad\qquad\fbox

\parindent=0mm \vspace{.4in}

{\bf Property 3.}{\it (Linearity) Let $\Psi$ be a non zero quaternion window function in $L^2(\mathbb R^2,\mathbb H)$ and $f_n\in L^2(\mathbb R^2,mathbb H),n\in\mathbb N$ then following holds: }

$$\mathbb S^{\mathbb H}_{\Psi,M_1,M_2}\left\{\sum_{n\in\mathbb N}\alpha_n f_n \right\}(\mathbf{u},\mathbf{w})=\sum_{n\in\mathbb N}\alpha_n\mathbb S^{\mathbb H}_{\Psi,M_1,M_2}[f_n](\mathbf{u},\mathbf{w}),\alpha_n\in\mathbb H\eqno(3.11)$$

\parindent=0mm \vspace{.1in}

{\bf Proof.} Follows directly from definition (3.2).\qquad\qquad\fbox

\parindent=0mm \vspace{.4in}

{\bf Property 4.}{\it(Parity) Let $\Psi\in L^2(\mathbb R^2,\mathbb H)$ be a non zero quaternion window function and $f\in L^2(\mathbb R^2,\mathbb H),$ then we have }
$$\mathbb S^{\mathbb H}_{P\Psi,M_1,M_2}[Pf](\mathbf{u},\mathbf{w})=\mathbb S^{\mathbb H}_{\Psi,M_1,M_2}[f](\mathbf{-u},\mathbf{-w})\eqno(3.12)$$
where $P\Psi(\mathbf{x})=\Psi(\mathbf{-x}).$

\parindent=0mm \vspace{.1in}

{\bf Proof.} Easy to prove so avoided.\qquad\qquad\fbox

\parindent=0mm \vspace{.4in}

{\bf Property 5.}{\it(Shift) Let $\Psi\in L^2(\mathbb R^2,\mathbb H)$ be a non zero quaternion window function and $f\in L^2(\mathbb R^2,\mathbb H),$ then following property holds }
$$\mathbb S^{\mathbb H}_{\Psi,M_1,M_2}[f(\mathbf{x-\alpha})](\mathbf{u},\mathbf{w})=e^{\frac{\mu_1}{2B_1}(A_1\alpha^2_1-2\alpha_1 w_1)}\mathbb S^{\mathbb H}_{\Psi,M_1,M_2}[\tilde{f}(\mathbf{u-\alpha},\mathbf{w})](\mathbf{.} ,\mathbf{.})e^{\frac{\mu_2}{2B_2}(A_2\alpha^2_2-2\alpha_2 w_2)}.\eqno(3.13)$$
where $\tilde{f}(\mathbf{t})=e^{{\mu_1\frac{A_1t_1\alpha_1}{B_1}}+{\mu_2\frac{A_2t_2\alpha_2}{B_2}}} f(\mathbf{t})$

\parindent=0mm \vspace{.1in}

{\bf Proof.}By taking into account of (3.2),we get

$$\mathbb S^{\mathbb H}_{\Psi,M_1,M_2}[f(\mathbf{x-\alpha})](\mathbf{u},\mathbf{w})=\int_{\mathbb R^2}K_{M_1}^{\mu_1}(x_1,w_1)f(\mathbf{x-\alpha})\overline{\Psi(\mathbf{u}-\mathbf{x},\mathbf{w})}K_{M_2}^{\mu_2}(x_2,w_2)d\mathbf{x}\eqno(3.14)$$
Applying change of variable $\mathbf{t}=\mathbf{x-\alpha}$ in (3.14),we obtain
\begin{align*}
\mathbb S^{\mathbb H}_{\Psi,M_1,M_2}[f(\mathbf{x}-\mathbf{\alpha})](\mathbf{u},\mathbf{w})&=\int_{\mathbb R^2}K_{M_1}^{\mu_1}(t_1+\alpha_1,w_1)f(\mathbf{t})\overline{\Psi(\mathbf{u}-(\mathbf{t+\alpha}),\mathbf{w})}K_{M_2}^{\mu_2}(t_2+\alpha_2,w_2)d\mathbf{t}\\\\
&=\int_{\mathbb R^2}\frac{1}{\sqrt{2\pi B_1}}e^{\mu_1\left(\frac{A_1}{2B_1}(t_1+\alpha_1)^2-\frac{(t_1+\alpha_1)w_1}{B_1}+\frac{D_1}{2B_1}w_1^2-\frac{\pi}{4} \right) }f(\mathbf{t})\overline{\Psi((\mathbf{u-\alpha})-\mathbf{t},\mathbf{w})}\\
&\qquad\qquad\times\frac{1}{\sqrt{2\pi B_2}}e^{\mu_2\left(\frac{A_2}{2B_2}(t_2+\alpha_2)^2-\frac{(t_2+\alpha_2)w_2}{B_2}+\frac{D_2}{2B_2}w_2^2-\frac{\pi}{4} \right) }d\mathbf{t}\\
&=\int_{\mathbb R^2}\frac{1}{\sqrt{2\pi B_1}}e^{\mu_1\left(\frac{A_1}{2B_1}(t^2_1+\alpha^2_1+2t_1\alpha_1)-\frac{(t_1+\alpha_1)w_1}{B_1}+\frac{D_1}{2B_1}w_1^2-\frac{\pi}{4} \right) }f(\mathbf{t})\overline{\Psi((\mathbf{u-\alpha})-\mathbf{t},\mathbf{w})}\\
&\qquad\qquad\times\frac{1}{\sqrt{2\pi B_2}}e^{\mu_2\left(\frac{A_2}{2B_2}(t^2_2+\alpha^2_2+t_2\alpha_2)-\frac{(t_2+\alpha_2)w_2}{B_2}+\frac{D_2}{2B_2}w_2^2-\frac{\pi}{4} \right) }d\mathbf{t}\\\\
&=e^{\frac{\mu_1}{2B_1}(A_1\alpha^2_1-2\alpha_1 w_1)}\int_{\mathbb R^2}\frac{1}{\sqrt{2\pi B_1}}e^{\mu_1\left(\frac{A_1}{2B_1}t^2_1-\frac{t_1w_1}{B_1}+\frac{D_1}{2B_1}w_1^2-\frac{\pi}{4} \right) }\\
&\qquad\qquad\times
\left\{e^{{\mu_1\frac{A_1t_1\alpha_1}{B_1}}+{\mu_2\frac{A_2t_2\alpha_2}{B_2}}} f(\mathbf{t})\right\}\overline{\Psi((\mathbf{u-\alpha})-\mathbf{t},\mathbf{w})}
\\
&\qquad\qquad\qquad\times\frac{1}{\sqrt{2\pi B_2}}e^{\mu_1\left(\frac{A_2}{2B_2}t^2_2-\frac{t_2w_2}{B_2}+\frac{D_2}{2B_2}w_2^2-\frac{\pi}{4} \right) }d\mathbf{t}e^{\frac{\mu_2}{B_2}(\frac{A_2\alpha^2_2}{2}-\alpha_2 w_2)}.
\end{align*}
This completes proof.\qquad\qquad\fbox

\parindent=0mm \vspace{.4in}

{\bf Property 6.}{\it(Modulation)Let $\Psi\in L^2(\mathbb R^2,\mathbb H)$ be a non zero quaternion window function and $f\in L^2(\mathbb R^2,\mathbb H),$ then we have }
\begin{align*}
\mathbb S^{\mathbb H}_{\Psi,M_1,M_2}[\mathcal M_\mathbf{s}f](\mathbf{u},\mathbf{w})&=e^{\mu_1\frac{D_1}{2}(2w_1s_1-B_1s^2_1)}\int_{\mathbb R^2}K_{M_1}^{\mu_1}(w_1-s_1B_1,x_1)f(\mathbf{x})\overline{\Psi((\mathbf{u-x}),\mathbf{w})}\\
&\qquad\qquad\qquad\qquad\qquad\times K_{M_2}^{\mu_2}(w_2-s_2B_2,x_2)d\mathbf{x}e^{\mu_2\frac{D_2}{2}(2w_2s_2-B_2s^2_2)},\tag{3.15}
\end{align*}

where $\mathcal M_{\mathbf{s}}f(\mathbf{x})=e^{\mu_1s_1x_1}f(\mathbf{x})e^{\mu_2s_2x_2}$\\

{\bf Proof.} It follows from(3.2) that,
\begin{align*}
\mathbb S^{\mathbb H}_{\Psi,M_1,M_2}[\mathcal M_\mathbf{s}f](\mathbf{u},\mathbf{w})&=\int_{\mathbb R^2}K_{M_1}^{\mu_1}(x_1,w_1)e^{\mu_1s_1x_1}f(\mathbf{x})e^{\mu_2s_2x_2}\overline{\Psi(\mathbf{u}-\mathbf{x},\mathbf{w})}K_{M_2}^{\mu_2}(x_2,w_2)d\mathbf{x}\\\\
&=\int_{\mathbb R^2}\frac{1}{\sqrt{2\pi B_1}}e^{\mu_1\left(\frac{A_1}{2B_1}x^2_1-\frac{x_1w_1}{B_1}+\frac{D_1}{2B_1}w_1^2+s_1x_1-\frac{\pi}{4} \right) }f(\mathbf{x})\overline{\Psi((\mathbf{u-x}),\mathbf{w})}\\
&\qquad\qquad\times \frac{1}{\sqrt{2\pi B_2}}e^{\mu_2\left(\frac{A_2}{2B_2}x^2_2-\frac{x_2w_2}{B_2}+\frac{D_2}{2B_2}w_2^2+s_2x_2-\frac{\pi}{4} \right) }d\mathbf{x}\\\\
&=\int_{\mathbb R^2}\frac{1}{\sqrt{2\pi B_1}}e^{\mu_1\left(\frac{A_1}{2B_1}x^2_1-\frac{x_1}{B_1}(w_1-s_1B_1)+\frac{D_1}{2B_1}w_1^2+-\frac{\pi}{4} \right) }f(\mathbf{x})\overline{\Psi((\mathbf{u-x}),\mathbf{w})}\\
&\qquad\qquad\times \frac{1}{\sqrt{2\pi B_2}}e^{\mu_2\left(\frac{A_2}{2B_2}x^2_2-\frac{x_2}{B_2}(w_2-s_2B_2)+\frac{D_2}{2B_2}w_2^2+s_2x_2-\frac{\pi}{4} \right) }d\mathbf{x}\\\\
&=\int_{\mathbb R^2}\frac{1}{\sqrt{2\pi B_1}}e^{\mu_1\left(\frac{A_1}{2B_1}x^2_1-\frac{x_1}{B_1}(w_1-s_1B_1)+\frac{D_1}{2B_1}((w_1-s_1B_1)+s_1B_1)^2-\frac{\pi}{4} \right) }f(\mathbf{x})\overline{\Psi((\mathbf{u-x}),\mathbf{w})}\\
&\qquad\qquad\times \frac{1}{\sqrt{2\pi B_2}}e^{\mu_2\left(\frac{A_2}{2B_2}x^2_2-\frac{x_2}{B_2}(w_2-s_2B_2)+\frac{D_2}{2B_2}((w_2-s_2B_2)+s_2B_2)^2-\frac{\pi}{4} \right) }d\mathbf{x}\\\\
&=\int_{\mathbb R^2}\frac{1}{\sqrt{2\pi B_1}}e^{\mu_1\left(\frac{A_1}{2B_1}x^2_1-\frac{x_1}{B_1}(w_1-s_1B_1)+\frac{D_1}{2B_1}((w_1-s_1B_1)^2-\frac{\pi}{4}+2(w_1-s_1B_1)s_1B_1+s_1B^2_1) \right) }\\
&\qquad\times f(\mathbf{x})\overline{\Psi((\mathbf{u-x}),\mathbf{w})}
\\
&\qquad\times \frac{1}{\sqrt{2\pi B_2}}e^{\mu_2\left(\frac{A_2}{2B_2}x^2_2-\frac{x_2}{B_2}(w_2-s_2B_2)+\frac{D_2}{2B_2}((w_2-s_2B_2)^2-\frac{\pi}{4}+2(w_2-s_2B_2)s_2B_2+s_2B^2_2) \right) } d\mathbf{x}\\\\
&=e^{\mu_1\frac{D_1}{2}(2w_1s_1-B_1s^2_1)}\int_{\mathbb R^2}\frac{1}{\sqrt{2\pi B_1}}e^{\mu_1\left(\frac{A_1}{2B_1}x^2_1-\frac{x_1}{B_1}(w_1-s_1B_1)+\frac{D_1}{2B_1}((w_1-s_1B_1)^2-\frac{\pi}{4}) \right) }\\
&\qquad\times f(\mathbf{x})\overline{\Psi((\mathbf{u-x}),\mathbf{w})}
\\
&\qquad\times \frac{1}{\sqrt{2\pi B_2}}e^{\mu_2\left(\frac{A_2}{2B_2}x^2_2-\frac{x_2}{B_2}(w_2-s_2B_2)+\frac{D_2}{2B_2}((w_2-s_2B_2)^2-\frac{\pi}{4}) \right) } d\mathbf{x}e^{\mu_2\frac{D_2}{2}(2w_2s_2-B_2s^2_2)}\\\\
&=e^{\mu_1\frac{D_1}{2}(2w_1s_1-B_1s^2_1)}\int_{\mathbb R^2}K_{M_1}^{\mu_1}(w_1-s_1B_1,x_1)f(\mathbf{x})\overline{\Psi((\mathbf{u-x}),\mathbf{w})}\\
&\qquad\qquad\qquad\qquad\qquad\times K_{M_2}^{\mu_2}(w_2-s_2B_2,x_2)d\mathbf{x}e^{\mu_2\frac{D_2}{2}(2w_2s_2-B_2s^2_2)}.
\end{align*}

Hence completes the proof.\qquad\qquad\fbox

\parindent=0mm \vspace{.1in}

Now we are going to develop the
orthogonality relation for Q-LCST and using this we derive
reconstruction formula associated with Q-LCST both of which are fundamental properties for signal analysis.Before we begin with orthogonality relation we note that throughout the rest of the paper, we will always assume that
  $$\int_{\mathbb R^2}|\Psi(\mathbf{u,w})|^2d\mathbf{u}=\Lambda_\Psi,\qquad0<\Lambda_\Psi<\infty\eqno(3.16)$$

\parindent=0mm \vspace{.1in}

{\bf Theorem 3.1(Orthogonality relation).}{\it Let $\Psi\in L^2(\mathbb R^2,\mathbb H)$ be a non zero quaternion window function and $f,g\in L^2(\mathbb R^2,\mathbb H),$ then we have }
$$\int_{\mathbb R^2}\int_{\mathbb R^2}\mathbb S^{\mathbb H}_{\Psi,M_1,M_2}[f](\mathbf{u},\mathbf{w})\overline{\mathbb S^{\mathbb H}_{\Psi,M_1,M_2}[g](\mathbf{u},\mathbf{w})}d\mathbf{w}d\mathbf{u}=\left\langle \Lambda_\Psi f,g\right\rangle\eqno(3.17)$$

\parindent=0mm \vspace{.1in}

{\bf Proof.} By definition of Q-LCST,we obtain
$$\begin{array}{lcr}
\displaystyle\int_{\mathbb R^2}\displaystyle\int_{\mathbb R^2}\mathbb S^{\mathbb H}_{\Psi,M_1,M_2}[f](\mathbf{u},\mathbf{w})\overline{\mathbb S^{\mathbb H}_{\Psi,M_1,M_2}[g](\mathbf{u},\mathbf{w})}d\mathbf{w}d\mathbf{u}&&\\\\
\qquad=\displaystyle\int_{\mathbb R^2}\displaystyle\int_{\mathbb R^2}\displaystyle\int_{\mathbb R^2}\frac{1}{\sqrt{2\pi B_1}}e^{-\mu_1\left(\frac{A_1}{2B_1}x^2_1-\frac{x_1w_1}{B_1}+\frac{D_1}{2B_1}w_1^2-\frac{\pi}{4} \right)}f(\mathbf{x})\overline{\Psi(\mathbf{u-x,w})}&&\\\\
\qquad\qquad\times \frac{1}{\sqrt{2\pi B_2}}e^{-\mu_2\left(\frac{A_2}{2B_2}x^2_2-\frac{x_2w_2}{B_2}+\frac{D_2}{2B_2}w_2^2-\frac{\pi}{4} \right)}
 \displaystyle\int_{\mathbb R^2}\frac{1}{\sqrt{2\pi B_1}}e^{\mu_1\left(\frac{A_1}{2B_1}t^2_1-\frac{t_1w_1}{B_1}+\frac{D_1}{2B_1}w_1^2-\frac{\pi}{4} \right)}\overline{g(\mathbf{t})}&&\\\\
 \qquad\qquad\times \Psi(\mathbf{u-t,w})\frac{1}{\sqrt{2\pi B_2}}e^{\mu_2\left(\frac{A_2}{2B_2}t^2_2-\frac{t_2w_2}{B_2}+\frac{D_2}{2B_2}w_2^2-\frac{\pi}{4} \right)}d\mathbf{t}d\mathbf{x}d\mathbf{u}d\mathbf{w}&&\\\\
 \qquad= \displaystyle\int_{\mathbb R^2}\displaystyle\int_{\mathbb R^2}\displaystyle\int_{\mathbb R^2}\displaystyle\int_{\mathbb R^2}\frac{1}{\sqrt{2\pi B_1}}e^{-\mu_1\frac{A_1}{2B_1}\left(x^2_1-t^2_1\right)}f(\mathbf{x})\overline{\Psi(\mathbf{u-x,w})}\frac{1}{\sqrt{2\pi B_2}}e^{-\mu_2\frac{A_2}{2B_2}\left(x^2_2-t^2_2\right)}&&\\\\
 \qquad\qquad\times\frac{1}{\sqrt{2\pi B_1}}e^{\frac{\mu_1}{B_1}\left(t_1-x_1\right)}\overline{g(\mathbf{t})}\Psi(\mathbf{u-x,w})\frac{1}{\sqrt{2\pi B_2}}e^{\frac{\mu_2}{B_2}\left(t_2-x_2\right)}d\mathbf{t}d\mathbf{x}d\mathbf{u}d\mathbf{w}&&\\\\
 \qquad= 2\pi B_1\displaystyle\int_{\mathbb R^2}\displaystyle\int_{\mathbb R^2}\displaystyle\int_{\mathbb R^2}\frac{1}{2\pi B_1}e^{-\mu_1\frac{A_1}{2B_1}\left(x^2_1-t^2_1\right)}f(\mathbf{x})\overline{\Psi(\mathbf{u-x,w})}\frac{1}{2\pi B_2}e^{-\mu_2\frac{A_2}{2B_2}\left(x^2_2-t^2_2\right)}&&\\\\
 \qquad\qquad\times\overline{g(\mathbf{t})}\Psi(\mathbf{u-t,w})\delta(\mathbf{t-x})d\mathbf{t}d\mathbf{x}d\mathbf{u}.2\pi B_2&&\\\\
 \qquad= \displaystyle\int_{\mathbb R^2}f(\mathbf{x})\overline{g(\mathbf{t})}d\mathbf{t}\displaystyle\int_{\mathbb R^2}\overline{\Psi(\mathbf{u-t,w})}\Psi(\mathbf{u-t,w})d\mathbf{u}&&\\\\
 \qquad= \displaystyle\int_{\mathbb R^2}f(\mathbf{x})\overline{g(\mathbf{t})}d\mathbf{t}\displaystyle\int_{\mathbb R^2}|{\Psi(\mathbf{u-t,w})}|^2d\mathbf{u}&&\\\\
 \qquad= \displaystyle\int_{\mathbb R^2}\Lambda_\Psi f(\mathbf{x})\overline{g(\mathbf{t})}d\mathbf{t}&&\\\\
  \qquad=\left\langle\Lambda_\Psi f,g\right\rangle.
\end{array}$$
Thus the proof is completed.\qquad\qquad\fbox

\parindent=0mm \vspace{.1in}

{\bf Remark 3.1.} {\it If we take $f=g$ in (3.17) Theorem 3.1 takes the form}
$$\int_{\mathbb R^2}\int_{\mathbb R^2}|\mathbb S^{\mathbb H}_{\Psi,M_1,M_2}[f](\mathbf{u},\mathbf{w})|^2d\mathbf{w}d\mathbf{u}=\Lambda_\Psi\|f\|^2.\eqno(3.18)$$

\parindent=0mm \vspace{.4in}

{\bf Theorem 3.2(Reconstruction formula).}{\it If $f\in L^2(\mathbb R^2,\mathbb H)$, then f can be reconstructed by the formula
 \begin{align*}f(\mathbf{x})&=\frac{1}{\Lambda_\Psi}\int_{\mathbb R^2}\int_{\mathbb R^2}\frac{1}{\sqrt{2\pi B_1}}e^{-\mu_1(\frac{A_1}{2B_1}x^2_1-\frac{x_1}{B_1}w_1+\frac{D_1}{2B_1}w^2_1-\frac{\pi}{4})}\mathbb S^{\mathbb H}_{\Psi,M_1,M_2}[f](\mathbf{u},\mathbf{w})\\
 &\qquad\qquad\qquad\times\Psi(\mathbf{u-x,w})\frac{1}{\sqrt{2\pi B_2}}e^{-\mu_2(\frac{A_2}{2B_2}x^2_2-\frac{x_2}{B_2}w_2+\frac{D_2}{2B_2}w^2_2-\frac{\pi}{4})}d\mathbf{w}d\mathbf{u}.\tag{3.19}
\end{align*}
Where $\Psi\in L^2(\mathbb R^2,\mathbb H)$ is a quaternion window function that satisfies (3.16)}

\parindent=0mm \vspace{.1in}

{\bf Proof.} From theorem {3.1}, we have

$$\begin{array}{lcr}
\left\langle f\Lambda_\Psi,g\right\rangle &&\\\\
\qquad=\displaystyle\int_{\mathbb R^2}\displaystyle\int_{\mathbb R^2}\mathbb S^{\mathbb H}_{\Psi,M_1,M_2}[f](\mathbf{u},\mathbf{w})\overline{\mathbb S^{\mathbb H}_{\Psi,M_1,M_2}[g](\mathbf{u},\mathbf{w})}d\mathbf{w}d\mathbf{u}&&\\\\
\qquad=\displaystyle\int_{\mathbb R^2}\displaystyle\int_{\mathbb R^2}\mathbb S^{\mathbb H}_{\Psi,M_1,M_2}[f](\mathbf{u},\mathbf{w})\displaystyle\int_{\mathbb R^2}e^{-\mu_1(\frac{A_1}{2B_1}x^2_1-\frac{x_1}{B_1}w_1+\frac{D_1}{2B_1}w^2_1-\frac{\pi}{4})}\overline{g(\mathbf{x})}&&\\
\qquad\qquad\qquad\qquad\times\Psi(\mathbf{u-x,w})\frac{1}{\sqrt{2\pi B_2}}e^{-\mu_2(\frac{A_2}{2B_2}x^2_2-\frac{x_2}{B_2}w_2+\frac{D_2}{2B_2}w^2_2-\frac{\pi}{4})}d\mathbf{w}d\mathbf{u}d\mathbf{x}&&\\\\
\qquad=\displaystyle\int_{\mathbb R^2}\displaystyle\int_{\mathbb R^2}\displaystyle\int_{\mathbb R^2}\mathbb S^{\mathbb H}_{\Psi,M_1,M_2}[f](\mathbf{u},\mathbf{w})e^{-\mu_1(\frac{A_1}{2B_1}x^2_1-\frac{x_1}{B_1}w_1+\frac{D_1}{2B_1}w^2_1-\frac{\pi}{4})}&&\\
\qquad\qquad\qquad\qquad\times\Psi(\mathbf{u-x,w})\frac{1}{\sqrt{2\pi B_2}}e^{-\mu_2(\frac{A_2}{2B_2}x^2_2-\frac{x_2}{B_2}w_2+\frac{D_2}{2B_2}w^2_2-\frac{\pi}{4})}\overline{g(\mathbf{x})}d\mathbf{w}d\mathbf{u}d\mathbf{x}&&\\\\
\qquad=\left\langle\displaystyle\int_{\mathbb R^2}\displaystyle\int_{\mathbb R^2}\mathbb S^{\mathbb H}_{\Psi,M_1,M_2}[f](\mathbf{u},\mathbf{w})e^{-\mu_1(\frac{A_1}{2B_1}x^2_1-\frac{x_1}{B_1}w_1+\frac{D_1}{2B_1}w^2_1-\frac{\pi}{4})}\right.&&\\
\qquad\qquad\qquad\times\left.\Psi(\mathbf{u-x,w})\frac{1}{\sqrt{2\pi B_2}}e^{-\mu_2(\frac{A_2}{2B_2}x^2_2-\frac{x_2}{B_2}w_2+\frac{D_2}{2B_2}w^2_2-\frac{\pi}{4})}d\mathbf{w}d\mathbf{u},g\right\rangle&&\\\\
\end{array}$$
As last equation is true for every $g\in L^2(\mathbb R^2,\mathbb H)$ therefore we have
\begin{align*}
f(\mathbf{x})\Lambda_\Psi&=\int_{\mathbb R^2}\int_{\mathbb R^2}\frac{1}{\sqrt{2\pi B_1}}e^{-\mu_1(\frac{A_1}{2B_1}x^2_1-\frac{x_1}{B_1}w_1+\frac{D_1}{2B_1}w^2_1-\frac{\pi}{4})}\mathbb S^{\mathbb H}_{\Psi,M_1,M_2}[f](\mathbf{u},\mathbf{w})\\
 &\qquad\qquad\qquad\times\Psi(\mathbf{u-x,w})\frac{1}{\sqrt{2\pi B_2}}e^{-\mu_2(\frac{A_2}{2B_2}x^2_2-\frac{x_2}{B_2}w_2+\frac{D_2}{2B_2}w^2_2-\frac{\pi}{4})}d\mathbf{w}d\mathbf{u}.
 \end{align*}
Which implies
\begin{align*}
 f(\mathbf{x})&=\frac{1}{\Lambda_\Psi}\int_{\mathbb R^2}\int_{\mathbb R^2}\frac{1}{\sqrt{2\pi B_1}}e^{-\mu_1(\frac{A_1}{2B_1}x^2_1-\frac{x_1}{B_1}w_1+\frac{D_1}{2B_1}w^2_1-\frac{\pi}{4})}\mathbb S^{\mathbb H}_{\Psi,M_1,M_2}[f](\mathbf{u},\mathbf{w})\\
 &\qquad\qquad\qquad\qquad\qquad\times\Psi(\mathbf{u-x,w})\frac{1}{\sqrt{2\pi B_2}}e^{-\mu_2(\frac{A_2}{2B_2}x^2_2-\frac{x_2}{B_2}w_2+\frac{D_2}{2B_2}w^2_2-\frac{\pi}{4})}d\mathbf{w}d\mathbf{u}.
 \end{align*}
 Hence completes the proof\qquad\qquad\fbox

\parindent=0mm \vspace{.4in}

 {\bf{4. Uncertainty Principles for the Q-LCST   }}

 \parindent=0mm \vspace{.1in}

The classical uncertainty principle of harmonic analysis states
that a non-trivial function and its Fourier transform cannot
both be simultaneously sharply localized. In quantum mechanics,
an uncertainty principle asserts that one cannot be certain of
the position and of the velocity of an electron (or any particle)
at the same time. In other words, increasing the knowledge of
the position decreases the knowledge of the velocity or momentum
of an electron.The uncertainty principles in harmonic analysis are of central importance as they provide
a lower bound for optimal simultaneous resolution in the time and frequency domains.In this Section, we shall establish an analogue of the well-known Heisenberg’s uncertainty
inequality and the corresponding logarithmic uncertainty principle for
the Q-LCST as defined by (3.2). First, we prove
the following lemma.

\parindent=0mm \vspace{.2in}

{\bf Lemma 4.1.} {\it Let $\Psi\in L^2(\mathbb R^2,\mathbb H)$ be a non zero  quaternion window function, then for every
$f \in L^2(\mathbb R^2,\mathbb H)$, we have
$$\Lambda_\Psi\int_{\mathbb R^2}x^2_s|f(\mathbf{x})|^2d\mathbf{x}=\int_{\mathbb R^2}\int_{\mathbb R^2}x^2_s|\mathcal L^{-1}_{M_1,M_2}\{\mathbb S^{\mathbb H}_{\Psi,M_1,M_2}[f](\mathbf{u},\mathbf{w})\}(\mathbf{x})|^2d\mathbf{x}d\mathbf{u}.\eqno(4.1)$$
where $s=1,2.$}

\parindent=0mm \vspace{.1in}

{\bf Proof.} We avoided proof as it follows by theorem (3.1) and theorem (3.2).\qquad\qquad\fbox

\parindent=0mm \vspace{.2in}

 Now we can establish the Heisenberg-type inequalities for the proposed Q-LCST as defined by (3.2).

 \parindent=0mm \vspace{.2in}

{\bf Theorem 4.2.} {\it Let  $\Psi\in L^2(\mathbb R^2,\mathbb H)$ be a non zero  quaternion window function and $\mathbb S^{\mathbb H}_{\Psi,M_1,M_2}[f]$ be the Q-LCST of any signal $f\in L^2(\mathbb R^2,\mathbb H),$ then we have
$$\left(\int_{\mathbb R^2}\int_{\mathbb R^2}w^2_s|\mathbb S^{\mathbb H}_{\Psi,M_1,M_2}[f](\mathbf{u,w})|^2d\mathbf{u}d\mathbf{w}\right)^{1/2}\left(\int_{\mathbb R^2}x^2_s|f(\mathbf{x})|^2d\mathbf{x}\right)^{1/2}\ge\frac{n_s\sqrt{\Lambda_\Psi}}{2}\|f\|^2\eqno(4.2)$$
where $s=1,2.$}

\parindent=0mm \vspace{.1in}

 {\bf Proof.} By virtue of the Heisenberg’s inequality for the QLCT
 \cite{14QWLCT}, we can write
 $$\int_{\mathbb R^2}x^2_s|f(\mathbf{x})|^2d\mathbf{x}\int_{\mathbb R^2}w^2_s|\mathcal L^{\mathbb H}_{M_1,M_2}[f](w)|^2d\mathbf{w}\ge\frac{n^2_s}{4}\left(\int_{\mathbb R^2}|f|^2d\mathbf{x}\right)^2.\eqno(4.3)$$
 Equation(4.3) can be rewritten as
 $$\int_{\mathbb R^2}x^2_s|\mathcal L^{-1}_{M_1,M_2}\{\mathcal L^{\mathbb H}_{M_1,M_2}[f]\}|^2d\mathbf{x}\int_{\mathbb R^2}w^2_s|\mathcal L^{\mathbb H}_{M_1,M_2}[f](\mathbf{w})|^2d\mathbf{w}\ge\frac{n^2_s}{4}\left(\int_{\mathbb R^2}|f|^2d\mathbf{x}\right)^2.\eqno(4.4)$$
 Now by applying Plancherel’s theorem for the QLCT
to the right-hand side of (4.4), we have
$$\int_{\mathbb R^2}x^2_s|\mathcal L^{-1}_{M_1,M_2}\{\mathcal L^{\mathbb H}_{M_1,M_2}[f]\}|^2d\mathbf{x}\int_{\mathbb R^2}w^2_s|\mathcal L^{\mathbb H}_{M_1,M_2}[f](\mathbf{w})|^2d\mathbf{w}\ge\left(\frac{n_s}{2}\int_{\mathbb R^2}|\mathcal L^{\mathbb H}_{M_1,M_2}[f](\mathbf{w})|^2d\mathbf{w}\right)^2.\eqno(4.5)$$
Since $\mathbb S^{\mathbb H}_{\Psi,M_1,M_2}[f]\in L^2(\mathbb R^2,\mathbb H)$,therefore we can replace $\mathcal L^{\mathbb H}_{M_1,M_2}[f]$ by  $\mathbb S^{\Psi,\mathbb H}_{\Psi,M_1,M_2}[f]$ on the both sides of (4.5) to get
$$\int_{\mathbb R^2}x^2_s|\mathcal L^{-1}_{M_1,M_2}\{\mathbb S^{\mathbb H}_{\Psi,M_1,M_2}[f](\mathbf{u,w})\}|^2d\mathbf{x}\int_{\mathbb R^2}w^2_s|\mathbb S^{\Psi,\mathbb H}_{M_1,M_2}[f](\mathbf{u,w})|^2d\mathbf{w}\ge\left(\frac{n_s}{2}\int_{\mathbb R^2}|\mathbb S^{\mathbb H}_{\Psi,M_1,M_2}[f]\mathbf{w}|^2d\mathbf{w}\right)^2.\eqno(4.6)$$
On taking square root to both sides of (4.6) and  integrating   with respect $d\mathbf{u}$,we have
\begin{align*}
\int_{\mathbb R^2}\left\{\left(\int_{\mathbb R^2}x^2_s|\mathcal L^{-1}_{M_1,M_2}\{\mathbb S^{\mathbb H}_{\Psi,M_1,M_2}[f](\mathbf{u,w})\}|^2d\mathbf{x}\right)^{1/2}\left(\int_{\mathbb R^2}w^2_s|\mathbb S^{\mathbb H}_{\Psi,M_1,M_2}[f](\mathbf{w})|^2d\mathbf{w}\right)^{1/2}\right\}d\mathbf{u}\\
\qquad\ge\frac{n_s}{2}\int_{\mathbb R^2}\int_{\mathbb R^2}|\mathbb S^{\mathbb H}_{\Psi,M_1,M_2}[f]\mathbf{w}|^2d\mathbf{w}d\mathbf{u}.\tag{4.7}
\end{align*}
Now by virtue of Cauchy-Schwarz inequality (4.7) becomes
\begin{align*}
\left(\int_{\mathbb R^2}\int_{\mathbb R^2}x^2_s|\mathcal L^{-1}_{M_1,M_2}\{\mathbb S^{\mathbb H}_{\Psi,M_1,M_2}[f](\mathbf{u,w})\}|^2d\mathbf{x}d\mathbf{u}\right)^{1/2}\left(\int_{\mathbb R^2}\int_{\mathbb R^2}w^2_s|\mathbb S^{\mathbb H}_{\Psi,M_1,M_2}[f](\mathbf{w})|^2d\mathbf{w}d\mathbf{u}\right)^{1/2}\\
\qquad\ge\frac{n_s}{2}\int_{\mathbb R^2}\int_{\mathbb R^2}|\mathbb S^{\mathbb H}_{\Psi,M_1,M_2}[f]\mathbf{w}|^2d\mathbf{w}d\mathbf{u}.\tag{4.8}
\end{align*}
Now by applying Lemma 4.1 on L.H.S and Remark 3.1 on R.H.S of the above
inequality,we have
\begin{align*}
\left(\Lambda_\Psi\int_{\mathbb R^2}x^2_s|f(\mathbf{x})|^2d\mathbf{x}\right)^{1/2}\left(\int_{\mathbb R^2}\int_{\mathbb R^2}w^2_s|\mathbb S^{\mathbb H}_{\Psi,M_1,M_2}[f](\mathbf{u,w})|^2d\mathbf{w}d\mathbf{u}\right)^{1/2}\\
\qquad\ge\frac{n_s\Lambda_\Psi}{2}\|f\|^2.\tag{4.9}
\end{align*}
 On further simplifying (4.9),we get
 $$\left(\Lambda_\Psi\int_{\mathbb R^2}x^2_s|f(\mathbf{x})|^2d\mathbf{x}\right)^{1/2}\left(\int_{\mathbb R^2}\int_{\mathbb R^2}w^2_s|\mathbb S^{\mathbb H}_{\Psi,M_1,M_2}[f](\mathbf{u,w})|^2d\mathbf{w}d\mathbf{u}\right)^{1/2}\ge\frac{n_s\Lambda_\Psi}{2}\|f\|^2.\eqno(4.10)$$
 Hence,
 $$\left(\int_{\mathbb R^2}\int_{\mathbb R^2}w^2_s|\mathbb S^{\mathbb H}_{\Psi,M_1,M_2}[f](\mathbf{u,w})|^2d\mathbf{w}d\mathbf{u}\right)^{1/2}\left(\int_{\mathbb R^2}x^2_s|f(\mathbf{x})|^2d\mathbf{x}\right)^{1/2}\ge\frac{n_s\sqrt{\Lambda_\Psi}}{2}\|f\|^2.\eqno(4.11)$$
 Which completes the proof.\qquad\qquad\fbox

 \parindent=0mm \vspace{.4in}

 We now establish the logarithmic uncertainty principle for the Q-LCST
 as defined by (3.2).

 \parindent=0mm \vspace{.1in}

 {\bf Theorem 4.3} {\it For  a quaternion window function  $\Psi\in\mathcal S(\mathbb R^2,\mathbb H)$  and for $f\in \mathcal S(\mathbb R^2,\mathbb H)$ the Q-LCST satisfies the following logarithmic estimate of the uncertainty inequality:
 \begin{align*}
 \int_{\mathbb R^2}\int_{\mathbb R^2}\ln|\mathbf{w}||\mathbb S^{\mathbb H}_{\Psi,M_1,M_2}[f](\mathbf{u,w})|^2d\mathbf{w}d\mathbf{u}+\Lambda_\Psi\int_{\mathbb R^2}\ln|\mathbf{x}||f(\mathbf{x})|^2d\mathbf{x}\ge\mathcal D\Lambda_\Psi\|f\|^2\tag{4.12}
 \end{align*}
 where $\mathcal D=\left(\frac{\Gamma'(1/2)}{\Gamma(1/2)}-\ln 2\right)$ and $\Gamma$ is a Gamma function.}

\parindent=0mm \vspace{.1in}

 {\bf Proof.}  By virtue of the logarithmic inequality for the QLCT
 \cite{14QWLCT}, we can write
 $$\int_{\mathbb R^2}\ln|\mathbf{w}||\mathcal L^{\mathbb H}_{M_1,M_2}[f](\mathbf{w})|^2d\mathbf{w}+\int_{\mathbb R^2}\ln|\mathbf{x}||f(\mathbf{x})|^2d\mathbf{x}\ge\mathcal D\int_{\mathbb R^2}|f(\mathbf{x})|^2d\mathbf{x}.\eqno(4.13)$$\\
 Now invoking the inversion formula of QLCT on the L.H.S and Parseval’s formula for QLCT
on R.H.S, we obtain
\begin{align*}
\int_{\mathbb R^2}\ln|\mathbf{x}||\mathcal L^{-1}_{M_1,M_2}\{\mathcal L^{\mathbb H}_{M_1,M_2}[f]\}(\mathbf{x})|^2d\mathbf{x}+\int_{\mathbb R^2}\ln|\mathbf{w}||\mathcal L^{\mathbb H}_{M_1,M_2}[f] (\mathbf{w})|^2d\mathbf{w}\\
\qquad\ge\mathcal D\int_{\mathbb R^2}|\mathcal L^{\mathbb H}_{M_1,M_2}[f](\mathbf{x})|^2d\mathbf{x}\tag{4.14}.\\
\end{align*}
Since $\mathcal L^{\mathbb H}_{M_1,M_2}[f]$ and $\mathbb S^{\mathbb H}_{\Psi ,M_1,M_2}[f]$ are in $\mathcal S(\mathbb R^2,\mathbb H)$ thus we can replace $\mathcal L^{\mathbb H}_{M_1,M_2}[f]$ by  $\mathbb S^{\Psi,\mathbb H}_{\Psi,M_1,M_2}[f]$ on the both sides of (4.14) to get

\begin{align*}
\int_{\mathbb R^2}\ln|\mathbf{w}||\mathbb S^{\mathbb H}_{\Psi,M_1,M_2}[f](\mathbf{u,w})|^2d\mathbf{w}+\int_{\mathbb R^2}\ln|\mathbf{x}||\mathcal L^{-1}_{M_1,M_2}\{\mathbb S^{\mathbb H}_{\Psi,M_1,M_2}[f](\mathbf{u,w})\}(\mathbf{x})|^2d\mathbf{x}\\
\qquad\ge\mathcal D\int_{\mathbb R^2}|\mathbb S^{\mathbb H}_{\Psi,M_1,M_2}[f](\mathbf{u,w})|^2d\mathbf{w}\tag{4.15}.\\
\end{align*}
On integrating (4.15) with respect du and then by applying Fubini theorem,we have
\begin{align*}
\int_{\mathbb R^2}\int_{\mathbb R^2}\ln|\mathbf{w}||\mathbb S^{\mathbb H}_{\Psi,M_1,M_2}[f](\mathbf{u,w})|^2d\mathbf{w}d\mathbf{u}+\int_{\mathbb R^2}\int_{\mathbb R^2}\ln|\mathbf{x}||\mathcal L^{-1}_{M_1,M_2}\{\mathbb S^{\mathbb H}_{\Psi,M_1,M_2}[f](\mathbf{u,w})\}(\mathbf{x})|^2d\mathbf{x}d\mathbf{u}\\
\qquad\ge\mathcal D\int_{\mathbb R^2}\int_{\mathbb R^2}|\mathbb S^{\mathbb H}_{\Psi,M_1,M_2}[f](\mathbf{u,w})|^2d\mathbf{w}d\mathbf{u}\tag{4.16}.\\
\end{align*}
Now applying Lemma 4.1 on L.H.S and Remark 3.1 on R.H.S of (4.16), we get
 $$
\int_{\mathbb R^2}\int_{\mathbb R^2}\ln|\mathbf{w}||\mathbb S^{\mathbb H}_{\Psi,M_1,M_2}[f](\mathbf{u,w})|^2d\mathbf{w}d\mathbf{u}+\Lambda_\Psi\int_{\mathbb R^2}\ln|\mathbf{x}|f(\mathbf{x})|^2d\mathbf{x}
\ge\mathcal D\Lambda_\Psi\|f\|^2.$$

Which completes the proof.\qquad\qquad\fbox

\parindent=0mm \vspace{.1in}

{\bf Remark 4.1} {\it If f is normalized, namely $\|f\|_{L^2(\mathbb R^2,\mathbb H)}=1$ then (4.12) becomes}
\begin{align*}
 \int_{\mathbb R^2}\int_{\mathbb R^2}\ln|\mathbf{w}||\mathbb S^{\mathbb H}_{\Psi,M_1,M_2}[f](\mathbf{u,w})|^2d\mathbf{w}d\mathbf{u}+\Lambda_\Psi\int_{\mathbb R^2}\ln|\mathbf{x}||f(\mathbf{x})|^2d\mathbf{x}\ge\mathcal D\Lambda_\Psi
 \end{align*}

\parindent=0mm \vspace{.2in}
{\bf{Conclusions}}

\parindent=0mm \vspace{.2in}
We have introduced the notion of quaternion linear canonical S-transform(Q-LCST) which is an extension of the linear canonical S-transform. Firstly, we discussed  the fundamental properties of quaternion linear canonical S-transform(Q-LCST) and then establish some basic results including orthogonality relation and reconstruction formula. Most importantly, we  derived  the stronger version of the associated Heisenberg’s uncertainty inequality and the corresponding logarithmic version for quaternion linear canonical S-transform(Q-LCST). This will pay way for further findings in this area. Our research will help in finding more generalised and stronger versions of the inequality which will revolutionise the signal and image processing.

\parindent=0mm \vspace{.2in}
{\bf{Declarations}}
\begin{itemize}
\item  Availability of data and materials: The data is provided on the request to the authors.
\item Competing interests: The authors have no competing interests.
\item Funding: No funding was received for this work
\item Author's contribution: Both the authors equally contributed towards this work.
\item Acknowledgements: This work  is supported by the UGC-BSR Research Start Up Grant(No. F.30-498/2019(BSR)) provided by UGC, Govt. of India.

\end{itemize}

\parindent=0mm \vspace{.2in}
{\bf{References}}

\end{document}